\documentclass[a4paper,twoside,10pt,reqno]{amsart}
\usepackage[a4paper,left=3cm,right=3cm, top=3cm, bottom=3cm]{geometry}
\usepackage[latin1]{inputenc}
\usepackage{mathrsfs}
\usepackage{enumerate}
\usepackage{graphicx}
\usepackage{epstopdf}
\usepackage{amsmath}
\usepackage{amsthm}
\usepackage{algorithm,algpseudocode}
\usepackage{amssymb}
\usepackage{tikz}
\usepackage{stmaryrd}
\usepackage{subfigure}
\usepackage{float}
\usepackage{bigints}
\usepackage{cite}
\usepackage{color}
\usepackage[abs]{overpic}
\usepackage[font=footnotesize,labelfont=bf]{caption}
\usepackage{cases}
\usepackage[author={Lorenzo}]{pdfcomment}
\usepackage{comment}
\usepackage{algorithm,algpseudocode}
\usepackage{xcolor}

\newcommand{\Vertiii}[1]{{\left\vert\kern-0.25ex\left\vert\kern-0.25ex\left\vert #1 
		\right\vert\kern-0.25ex\right\vert\kern-0.25ex\right\vert}}

\restylefloat{table}
\theoremstyle{plain}

\newcommand{\h}{h}
\newcommand{\p}{p}
\newcommand{\taun}{\mathcal T_n}
\newcommand{\Ical}{\mathcal I}
\newcommand{\E}{K}
\newcommand{\Ehat}{\widehat \E}

\newcommand{\F}{F}
\newcommand{\Pbb}{\mathbb P}
\newcommand{\Rbb}{\mathbb R}
\newcommand{\Nbb}{\mathbb N}
\newcommand{\bbold}{\mathbf b}
\newcommand{\bboldz}{\bbold_0}
\renewcommand{\c}{c}
\newcommand{\cz}{\c_0}
\newcommand{\cs}{\c_s}
\newcommand{\f}{f}
\newcommand{\gD}{g_D}

\newcommand{\Vn}{V_n}
\newcommand{\un}{u_n}
\newcommand{\vn}{v_n}
\newcommand{\Bn}{B_n}

\newcommand{\Qbb}{\mathbb Q}
\newcommand{\PhiE}{\Phi_\E}

\newcommand{\Gammain}{\Gamma_{-}}
\newcommand{\Gammaout}{\Gamma_{+}}

\newcommand{\nbf}{\mathbf n}
\newcommand{\nbfGamma}{\nbf_\Gamma}
\newcommand{\hE}{\h_\E}
\newcommand{\nbfE}{\nbf_\E}
\newcommand{\GammaEin}{\Gammain^\E}
\newcommand{\GammaEout}{\Gammaout^\E}
\newcommand{\GammaE}{\Gamma^\E}
\newcommand{\Fn}{F_n}
\newcommand{\Pip}{\Pi_\p}
\newcommand{\DG}{\text{dG}}
\newcommand{\bo}{b_1}
\newcommand{\bt}{b_2}
\newcommand{\bd}{b_d}

\newcommand{\ujump}[1]{\lfloor #1\rfloor}   

\author{Zhaonan Dong \and Lorenzo Mascotto}
\address[Z. Dong]{ Inria, 2 rue Simone Iff, 75589 Paris, France and
	CERMICS, Ecole des Ponts, 77455 Marne-la-Vall\'{e}e 2, France}
\email{zhaonan.dong@inria.fr}
\address[L. Mascotto]{Fakult\"at f\"ur Mathematik, Universit\"at Wien, 1090 Vienna, Austria}
\email{lorenzo.mascotto@univie.ac.at}
\date{}
\title[On the suboptimality of the $p$-version dG methods for hyperbolic problems]{On the suboptimality of the $p$-version discontinuous Galerkin methods for first order hyperbolic problems}

\begin{document}
%%%%%%%%%%%%%%%%%%%%%%%%%%%%%%%%%%%%%
	
\begin{abstract}
\noindent
We address the issue of the suboptimality in the $p$-version discontinuous Galerkin (dG) methods for first order hyperbolic problems.
The convergence rate is derived for the upwind dG scheme on tensor product meshes in any dimension.
The standard proof in seminal work~\cite{HSS_hpDG} leads to suboptimal convergence in terms of the polynomial degree by~$3/2$ order for general convection fields, with the exception of piecewise multi-linear convection fields, which rather yield optimal convergence.
Such suboptimality is not observed numerically.
Thus, it might be caused by a limitation of the analysis, which we partially overcome:
for a special class of convection fields,
we shall show that the dG method has a $p$-convergence rate suboptimal by~$1/2$ order only.
\medskip

\begin{center}
\textbf{SUBMITTED AS A PROCEEDING ARTICLE}
\end{center}

\medskip

\medskip
		
\end{abstract}
	
\maketitle
%----------------------------
\section{Introduction}
%----------------------------

Discontinuous Galerkin (dG) finite element methods were introduced in the early 1970s for the numerical solution of first-order hyperbolic problems \cite{reedhill} and the weak imposition of inhomogeneous boundary conditions for elliptic problems \cite{nits}. In the past several decades, dG methods have enjoyed considerable success as a standard variational framework for the numerical solution of many classes of problems involving partial differential equations (PDEs); see, e.g., monographs~\cite{MR1842161,DiPietroErn,MR3363720} for reviews of some of the main developments of dG methods.
The interest in dG methods can be attributed to a number of factors, including the great flexibility in dealing with $hp$-adaptivity and general shaped elements~\cite{cangianigeorgoulishouston_hpDGFEM_polygon,DGpoly2,DGpolybook},
as well as in solving convection-dominated PDEs; see, e.g., early works~\cite{MR1399638,MR1702201} concerning hyperbolic conservation laws and convection-diffusion problems.

Due to missing tools in the analysis, the convergence rate always contains suboptimality in terms of the polynomial degree~$p$.
In~\cite{HSS_hyp}, the first optimal convergence rate of $hp$-dG methods is derived for linear convection problems by using the SUPG stabilisation.
However, the authors provide numerical evidence that the $hp$-optimal convergence rate is achieved even without such stabilisation.
In seminal work~\cite{HSS_hpDG}, based on (back then) novel optimal approximation results for the $L^2$-orthogonal projection, the $hp$-optimal convergence rate is derived for dG methods applied to hyperbolic problems, under the technical assumption that the convection field is piecewise linear.
Moreover, whenever the above assumption is violated, the theoretical analysis in~\cite{HSS_hpDG} leads to error bounds that are suboptimal in terms of~$p$ by $3/2$ order.
Such suboptimality is yet not observed in the numerical experiments. 
Over the last two decades, the above mentioned technical assumption became standard in $hp$-dG methods for convection-diffusion-reaction and hyperbolic problems; see, e.g., \cite{MR2377441,DGpoly2,DGpolybook,DGease}.
It is still an open question, whether the $p$-suboptimality for dG methods by~$3/2$ order is true or not in general. 

Our contribution represents a further step in shedding light on this issue. Notably, we present the a priori error analysis for $hp$-dG methods applied to pure hyperbolic problems employing a class of convection field, including nonpolynomial cases.
The new error is $h$-optimal and $p$-suboptimal by~$1/2$ order only.

The rest of the paper is organized as follows:
the continuous problem and its dG discretization are addressed in Section~\ref{section:continuous};
the classical analysis from~\cite{HSS_hyp} is re-elaborated in Section~\ref{section:standard-analysis},
whereas the improved bounds under suitable assumptions on the convection field are the topic of Section~\ref{section:special-analysis};
we collect the conclusions in Section~\ref{section:conclusion}.

Throughout, we employ a standard notation for Sobolev spaces~\cite{adamsfournier}.

%----------------------------
\section{The continuous problem and its dG formulation} \label{section:continuous}
%----------------------------

%%%%%%
\subsection{The continuous problem}
%%%%%%
	
Let $\Omega$ be a bounded polyhedral domain in~$\Rbb^d$, $d \in \Nbb$, with boundary~$\Gamma$.
We denote the unit outward normal vector to $\Gamma$ at $x\in \Gamma$ by~$\nbfGamma(x)$ and introduce the Fichera function~$\bbold  \cdot \nbf$ on~$\Gamma$ to define
\begin{equation} \label{splitting boundary}
\Gammain := \{x\in \Gamma: \bbold(x)  \cdot \nbf(x) < 0\} \quad \text{and}  \quad \Gammaout :=  \{x\in \Gamma: \bbold(x)  \cdot \nbf(x) \geq 0\}.
\end{equation} 
In the following, the sets~$\Gammain$ and~$\Gammaout$ are referred to as the  \emph{inflow} and \emph{outflow} boundary, respectively,
and clearly form a nonoverlapping partition of~$\Gamma$.

Given~$\bbold \in [W^{1,\infty}(\Omega)]^d$, $\c\in L^\infty(\Omega)$, $\f \in L^2(\Omega)$, and~$\gD \in L^{2}(\Gammain)$,
we consider the convection-reaction problem
\begin{equation} \label{strong-formulation}
\begin{aligned}
\bbold \cdot \nabla u + c\, u &= \f \quad \quad \!\! \text{in } \Omega,  \\
u&= \gD \quad \text{on } \Gammain.
\end{aligned}
\end{equation} 
Assuming the existence of a positive constant~$\cs$ satisfying
\begin{equation} \label{standard-assumption:well-posedness}
\cz^2 := \c - \frac12 \nabla \cdot \bbold \ge \cs \qquad \text{for a.e. } x\in \Omega,
\end{equation}
the well-posedness of problem~\eqref{strong-formulation} follows, e.g., as in~\cite{MR1825801}.
	
%%%%%%
\subsection{The dG formulation on quadrilateral/hexahedral meshes}
%%%%%%
We are interested in discretizing solutions to~\eqref{strong-formulation} by means of a dG finite element method.
To the aim, consider sequences of meshes~$\{ \taun \}_n$ consisting of tensor product elements, which can be defined through an affine mapping~$\PhiE$ on the reference $d$-dimensional cube element~$\Ehat:=(-1,1)^d$.
For sake of simplicity, we assume the all elements~$\E \in \taun$ are shape regular.
We fix a uniform polynomial degree~$\p \in \Nbb$ and denote the space of tensor polynomials of degree~$\p$ over~$\Ehat$ by~$\Qbb_\p(\Ehat)$.
Next, we introduce the dG space
\[
\Vn:= \{ \vn \in L^2(\Omega) \mid \vn{}_{|\E} \circ \PhiE \in \Qbb_p(\E) \; \forall \E \in \taun  \}.
\]
Given an element~$\E \in \taun$, we set its diameter and outward pointing normal by~$\hE$ and~$\nbfE$, respectively.
Given boundary~$\GammaE$ of element~$\E$, we split it into the \emph{inflow} and \emph{outflow} parts~$\GammaEin$ and~$\GammaEout$ defined as
\begin{equation} \label{splitting element boundary}
\GammaEin := \{x\in \GammaE: \bbold(x)  \cdot \nbfE(x) < 0\} \quad \text{and}  \quad  \GammaEout := \{x\in \GammaE: \bbold(x)  \cdot \nbfE(x) \geq 0\}.
\end{equation}	
Next, we define the classical \emph{upwind jump} operator.
Given an internal face~$\F$, let~$\E_1$ and~$\E_2$ be two elements in~$\taun$ sharing~$\F$. 
Without loss of generality, we assume that~$\E_1$ is such that~$\bbold \cdot \nbfE(x)<0$ for almost all~$x$ in~$\F$.
Then, we set
\begin{equation} \label{jump}
\ujump{v}_\F := (v{}_{|\E_1} - v{}_{|\E_2})_{|\F} = v{}^{+}- v{}^{-}
\quad \quad \forall v\in H^1(\Omega, \taun).
\end{equation}
In the rest of this work, when no confusion occurs, we shall write~$\nbf$ instead of~$\nbfE(x)$.
	
Following seminal work~\cite{HSS_hpDG}, we consider the upwind dG variational formulation of~\eqref{strong-formulation}.
More precisely, introduce the dG bilinear form
\begin{equation} \label{dG-bf}
\begin{split}
\Bn(\un,\vn)
&:= \sum_{\E \in \taun} \Big((\bbold \cdot \nabla \un, \vn)_{0,\E} + (\c\ \un,\vn)_{0,\E} \\
& \quad \quad \quad \quad  - ((\bbold \cdot \nbf) \ujump{\un}, \vn^+)_{\GammaEin  \setminus \Gammain}- ((\bbold \cdot \nbf) \un^+, \vn^+)_{\GammaEin  \cap \Gammain} \Big)
\end{split}
\end{equation}
and the discrete right-hand side
\begin{equation} \label{dG-RHS}
\Fn(\vn) := \sum_{\E \in \taun} \Big( (\f,\vn)_{0,\E} - ( (\bbold\cdot \nbf) \gD,\vn^+)_{\GammaEin \cap \Gammain} \Big) \quad \forall \vn \in \Vn.
\end{equation}
The dG method we consider reads
\begin{equation} \label{dG}
\text{find } \un \in \Vn \quad \text{such that} \quad \Bn(\un,\vn) = \Fn(\vn) \quad \vn \in \Vn.
\end{equation}
Furthermore, we introduce the dG norm
\begin{equation} \label{dG:norm}
\begin{split}
\Vertiii{\vn}_{\DG}^2
:=& \sum_{\E \in \taun} \Big( \Vert \cz \vn \Vert^2_{0,\E} +\frac12 \Vert \sqrt{|\bbold \cdot \nbf |} \vn^+ \Vert_{0,\GammaEin \cap \Gamma}^2
+\frac12 \Vert \sqrt{|\bbold \cdot \nbf |} \vn^+ \Vert_{0,\GammaEout \cap \Gamma}^2 \\
&  \quad \quad \quad +\frac12 \Vert \sqrt{|\bbold \cdot \nbf |} \ujump{\vn}  \Vert^2_{\GammaEin\setminus \Gammain} \Big).
\end{split}
\end{equation}
with $c_0$ defined as in \eqref{standard-assumption:well-posedness}.

It is easy to check that
\begin{equation} \label{coercivity}
\begin{split}
\Bn(\vn, \vn) = \Vertiii{\vn}_{\DG}^2 \quad\quad \forall \vn \in \Vn.
\end{split}
\end{equation}
The well-posedness of method~\eqref{dG} can be found, e.g., in~\cite[section 2.3]{DiPietroErn}.

%----------------------------
\section{Standard analysis and suboptimality in terms of the polynomial degree} \label{section:standard-analysis}
%----------------------------
In this section, we recall the convergence analysis for method~\eqref{dG} from \cite{HSS_hpDG} and where suboptimal estimates in terms of~$\p$ appear.
	
Preliminary, for all~$\E \in \taun$, introduce the $L^2$ projector~$\Pip:L^2(\E) \rightarrow \Qbb_\p(\E)$ through the affine mapping
and recall the standard $\h\p$-approximation estimates, see, e.g., \cite{HSS_hpDG,MR3936936}:
for any function $u\in H^{\ell}(K)$ on the given element~$\E\in \taun$, the follow relations hold
\begin{equation} \label{hp-L2-projection}
\begin{split}
&	\Vert u - \Pip u \Vert_{0,\E}\lesssim  \Big(\frac{\hE}{\p+1}\Big)^{s} \vert u \vert_{s,\E}, \\
&	\Vert u - \Pip u \Vert_{0,\GammaE}\lesssim  \Big(\frac{\hE}{\p+1}\Big)^{s-\frac{1}{2}} \vert u \vert_{s,\E},
\end{split}
\end{equation}
with $s:= \min\{p+1,\ell \}$.

It is easy to check that method~\eqref{dG} is consistent, whence the following Galerkin orthogonality follows:
\begin{equation} \label{Galerkin:orthogonality}
\Bn(u-\un, \vn) = 0 \quad\quad \forall \vn \in \Vn.
\end{equation}
Then, we split error~$u-\un$ into~$\eta + \xi$, where
\begin{equation} \label{error:splitting}
\eta:= u-\Pip u, \quad \quad \xi:= \Pip u - \un.
\end{equation}
Using Galerkin orthogonality~\eqref{Galerkin:orthogonality} and the properties of orthogonal projector~$\Pip$, we readily have the error equation
\begin{equation} \label{error-equation}
0=\Bn(u-\un, \xi)= \Bn(\eta, \xi) + \Bn(\xi, \xi) \quad \quad \Longrightarrow \quad \quad \Vertiii{\xi}^2_{\DG} = - \Bn(\eta, \xi).
\end{equation}
Since estimates on term~$\eta$ are standard, error equation~\eqref{error-equation} allows us to show a bound on term~$\xi$, on which we now focus on.
	
We begin by computing the following error splitting:
\begin{equation} \label{error-splitting}
\begin{split}
\Bn(\eta,\xi) 
& = \sum_{\E \in \taun} \Big( (\bbold\cdot \nabla \eta,\xi)_{0,\E} + (\c \eta, \xi)_{0,\E} \\
&\quad \quad \quad - ((\bbold\cdot \nbf) \ujump{\eta}, \xi^+)_{0,\GammaEin \setminus \Gammain} - ((\bbold\cdot \nbf)\eta^+, \xi^+)_{0,\GammaEin \cap \Gammain} \Big)\\
& = \sum_{\E \in \taun} \Big( ( (\c-\nabla \cdot \bbold) \eta,\xi)_{0,\E} - (\bbold\cdot \nabla \xi, \eta)_{0,\E} \\
& \quad \quad \quad+ ((\bbold\cdot \nbf) \ujump{\xi}, \eta^-)_{0,\GammaEin \setminus \Gammain} + ((\bbold\cdot \nbf)\xi^+, \eta^+)_{0,\GammaEout \cap \Gammaout} \Big) =: T_1+T_2 + T_3 +T_4.
\end{split}
\end{equation}
We show upper bounds for the four terms on the right-hand side of~\eqref{error-splitting} and anticipate that our analysis on term~$T_2$ will lead to suboptimal bounds in terms of~$\p$.
Under further assumptions on vector~$\bbold$, we shall exhibit improved $\p$-bounds in Section~\ref{section:special-analysis} below.
	
We begin with term~$T_1$. Using that~$\bbold \in [W^{1,\infty}(\Omega)]^2$, ~$\c \in L^{\infty}(\Omega)$, and assumption \eqref{standard-assumption:well-posedness}, we obtain
\begin{equation} \label{T1}
T_1     \le \sum_{\E \in \taun} \Vert \c - \nabla \cdot \bbold \Vert_{\infty,\E} \Vert \eta\Vert_{0,\E} \Vert \xi \Vert_{0,\E} 
\lesssim \Vert \eta \Vert_{0,\Omega} \Vert c_0 \xi \Vert_{0,\Omega} \lesssim \Vert \eta \Vert_{0,\Omega}\Vertiii{\xi}_{\DG}.
\end{equation}
As for terms~$T_3$ and~$T_4$, we have
\begin{equation} \label{T3+T4}
T_3+T_4 \lesssim (\sum_{\E\in \taun} \Vert (\bbold \cdot \nbfE)^{\frac12} \eta \Vert_{0,\GammaE}^2 )^{\frac12} \Vertiii{\xi}_{\DG}
\lesssim (\sum_{\E\in \taun} \Vert \eta \Vert_{0,\GammaE}^2 )^{\frac12} \Vertiii{\xi}_{\DG}.
\end{equation}
As for term~$T_2$, if we assume that~$\bbold \cdot \nabla \xi \in \Vn$  for all~$\xi \in \Vn$, then~$T_2=0$.
Thus, inserting~\eqref{T1} and~\eqref{T3+T4} in~\eqref{error-splitting} yields
\begin{equation} \label{estimate-xi:T2-zero}
\Vertiii{\xi}_{\DG} \lesssim (\sum_{\E \in \taun} (\Vert \eta \Vert^2_{0,\E} + \Vert \eta \Vert^2_{0,\GammaE})   )^{\frac12}
\lesssim \Big(\frac{\hE}{\p+1}\Big)^{s-\frac{1}{2}} \vert u \vert_{s,\Omega},
\end{equation}
Using a triangle inequality, and combining~\eqref{hp-L2-projection} with~\eqref{estimate-xi:T2-zero} leads to a $p$-optimal error estimate.
	
Next, we focus on the case of nonzero~$T_2$ to investigate the $p$-suboptimality.
Using the definition of~$\eta$ in~\eqref{error:splitting}, and notably the property of orthogonal projection~$\Pip$, we can write
\[
T_2     := \sum_{\E\in\taun} \int_\E (\bbold \cdot \nabla\xi) \eta =  \sum_{\E\in\taun} \int_\E (\bbold \cdot \nabla\xi - \bboldz \cdot \nabla\xi) \eta,
\]
where~$\bboldz$ is the vector average over every~$\E$ of~$\bbold$. We deduce
\[
T_2 \le \sum_{\E\in\taun} \Vert \bbold-\bboldz \Vert_{\infty,\E} \vert \xi \vert_{1,\E} \Vert \eta \Vert_{0,\E}.
\]
On each element~$\E \in \taun$, we have the following approximation property and $\h\p$-polynomial inverse inequality:
\[
\Vert \bbold-\bboldz \Vert_{\infty,\E} \lesssim \hE \vert \bbold \vert_{W^{1,\infty}(\E)},\quad \quad \vert \xi \vert_{1,\E} \lesssim \frac{\p^2}{\hE} \Vert \xi \Vert_{0,\E}.
\]
In the light of the two above bounds and~\eqref{standard-assumption:well-posedness}, we have the following bound on term~$T_2$:
\begin{equation} \label{T2-general}
T_2 \lesssim \p^2 \Vert \eta \Vert_{0,\Omega} \Vertiii{ \xi}_{\DG}.
\end{equation}
Inserting~\eqref{T1}, \eqref{T2-general}, and~\eqref{T3+T4} in~\eqref{error-splitting},
and using~\eqref{hp-L2-projection} yield
\begin{equation} \label{estimate-xi:T2-general}
\Vertiii{\xi}_{\DG}
\lesssim (1+\p^2) \Vert \eta \Vert_{0,\Omega} + (\sum_{\E\in\taun} \Vert \eta \Vert^2_{0,\GammaE})^{\frac12}
\lesssim \left(\frac{\h^s}{(\p+1)^{s-2}}+ \frac{\h^{s-\frac12}}{(\p+1)^{s-\frac12}}\right) \vert u \vert_{s,\Omega}.
\end{equation}
Using a triangle inequality, and combining~\eqref{hp-L2-projection} with~\eqref{estimate-xi:T2-general} leads to the following $\p$-suboptimal error estimate:
\begin{equation} \label{error-estimate-general}
\Vertiii{u-\un}_{\DG} \lesssim  \left(\frac{\h^{s-\frac12}}{(\p+1)^{s-2}}\right) \vert u \vert_{s,\Omega}.
\end{equation}
The above error bound is optimal in $h$ but suboptimal in terms of~$p$ by~$3/2$ order, which is in accordance with \cite[Remark~$3.13$]{HSS_hpDG}.
Notwithstanding, such suboptimality is not observed in practice; see, e.g., \cite[Numerical~Example~$1$]{HSS_hpDG}.

This motivates Section~\ref{section:special-analysis}, where we shall exhibit improved estimates in terms of~$\p$, under further assumptions on convection field~$\bbold$.

%----------------------------
\section{Improved bounds for special convection field} \label{section:special-analysis}
%----------------------------
In this section, we show improved $\p$-error estimates, under the following assumption on convection field~$\bbold$:
\begin{equation} \label{assumption:bbold}
		\bbold \in [W^{2,\infty} (\E)]^d \quad  \forall  \E\in \taun, \quad \quad
		\bbold = [\bo(x_1), \bt(x_2),\dots, \bd(x_d)]^T.
\end{equation}
Since the $j$-th, $j=1,\dots,d$, component of~$\bbold$ is assumed to be single-valued in the $x_j$ variable, without loss of generality, we can assume that element~$\E$ is the Cartesian product of intervals~$I_j:=(\alpha_j,\beta_j)$, $j=1,2,\dots,d$.

We can re-write term~$T_2$ in~\eqref{error-splitting} as
\[
T_2 = \sum_{j=1}^d \int_\E b_j (x_j) \partial_j \xi \ \eta. 
\]
Fix~$j=1,\dots,d$.
The $j$-th partial derivative of~$\xi$ is a tensor polynomial of degree~$\p-1$ along direction~$x_j$ and~$\p$ along the others.
Define~$\Ical_j b_j$ as the linear interpolant of~$b_j$ at the end-points of interval~$I_j$:
\begin{equation} \label{linear interpolant}
\big( b_j-\Ical_j b_j \Big)(\alpha_j) = 		\big( b_j-\Ical_j b_j \Big)(\beta_j) = 0.
\end{equation}
Then, we clearly have that~$\Ical_j b_j \ \partial_i \xi \in \Pbb_\p(I_j)$ and $\Ical_j b_j \ \partial_i \xi \in \Vn$.
Consequently, the definition of orthogonal projection~$\Pip$ allows us to write
\[
\int_\E b_j \partial_j \xi \ \eta = \int_\E (b_j - \Ical_j b_j)\partial_j \xi \ \eta.
\]
As~$\alpha_j$ and~$\beta_j$ denote the endpoints of interval~$I_j$, standard properties of one dimensional linear interpolation operators guarantee the existence of~$\widetilde x_j$ such that
\[
b_j (x_j) - \Ical_j b_j (x_j) = \frac{b_j^{(2)} (\widetilde x_j)}{2} (x-a_j)(x-b_j) \quad \quad \forall x_j \in I_j.
\]
Defining the standard quadratic bubble function on~$I_j$
\[
w_j := -(x-a_j)(x-b_j),
\]
we can thus write
\[
\int_\E (b_j - \Ical_j b_j)\partial_j \xi \ \eta = \int_\E \Big( \frac{b_j - \Ical_j b_j}{\sqrt{w_j}}	\Big) \big(\sqrt{w_j}\partial_j \xi \big)\ \eta
	\le 	\left\Vert \frac{b_j - \Ical_j b_j}{\sqrt{w_j} } \right\Vert_{\infty,\E} \Vert \sqrt{w_j} \partial_j \xi \Vert_{0,\E} \Vert \eta \Vert_{0,\E}.
\]
Thanks to assumption~\eqref{assumption:bbold} and an $\h\p$-polynomial inverse estimate involving bubbles,
see, e.g., \cite[Lemma 3.42]{verfurth2013posteriori-book}, we deduce
\[
\left\Vert \frac{b_j - \Ical_j b_j}{\sqrt{w_j}} \right\Vert_{\infty,\E} \lesssim \hE \vert b_j \vert_{W^{2,\infty}(\E)}, \quad \quad \Vert \sqrt{w_j} \partial_j \xi \Vert_{0,\E} \lesssim \sqrt{\p (\p+1)} \Vert \xi \Vert_{0,\E}.
\]
Collecting all the above estimates yields the following improved bound on term~$T_2$:
\begin{equation} \label{T2-special}
T_2 \lesssim \h\p \Vert \eta \Vert_{0,\Omega} \Vertiii{\xi}_{\DG}.
\end{equation}
Inserting~\eqref{T1}, \eqref{T2-special}, and~\eqref{T3+T4} in~\eqref{error-splitting} gives the bound
\[
\Vertiii{\xi}_{\DG} \lesssim \frac{\h^{s-\frac{1}{2}}}{(\p+1)^{s-1}} \vert u \vert_{s,\Omega}, 
\]
which, combined with a triangle inequality and~\eqref{hp-L2-projection}, eventually entails the error estimate
\[
\Vertiii{u-\un}_{\DG} \lesssim \frac{\h^{s-\frac{1}{2}}}{(\p+1)^{s-1}} \vert u \vert_{s,\Omega}.
\]
Compared with error estimate~\eqref{error-estimate-general}, the $\p$-suboptimality improved by one order.  

%%%%%%%%%%%%%%%%%%	
\section{Conclusion} \label{section:conclusion}
%%%%%%%%%%%%%%%%%%
Employing a special class of convection fields,
we derived an improved $hp$-bound for dG methods discretising linear hyperbolic problems.
The new error bound is suboptimal by~$1/2$ order in~$p$ only, which improves the $3/2$ suboptimal order presented in~\cite{HSS_hpDG}.
Needless to say, the results in this work do not provide full answers to the open questions in~\cite{HSS_hpDG}, notably on the mismatch between the theoretical and numerical results on $\p$-convergence.
However, they shed additional light on such issues and shows the possibility of deriving sharper $hp$-error bounds for dG methods.

%----------------------------
{\footnotesize
\bibliography{bibliogr}
}
\bibliographystyle{plain}
%----------------------------

\end{document}